\documentclass[12pt,a4paper]{article}
\usepackage{amssymb}
\usepackage{amsmath}
\usepackage{amsthm}
\usepackage[arrow,matrix]{xy}

\let\top=\undefined

\DeclareMathOperator{\codim}{codim}
\DeclareMathOperator{\dimv}{\bf dim}

\DeclareMathOperator{\Ext}{Ext}
\DeclareMathOperator{\GL}{GL}
\DeclareMathOperator{\Hom}{Hom}

\DeclareMathOperator{\ind}{ind}

\DeclareMathOperator{\rad}{rad}

\DeclareMathOperator{\Reg}{Reg}
\DeclareMathOperator{\rep}{rep}
\DeclareMathOperator{\soc}{soc}
\DeclareMathOperator{\Sing}{Sing}
\DeclareMathOperator{\top}{top}

\newcommand{\BA}{{\mathbb A}}

\newcommand{\BC}{{\mathbb C}}
\newcommand{\BD}{{\mathbb D}}
\newcommand{\BE}{{\mathbb E}}
\newcommand{\BM}{{\mathbb M}}
\newcommand{\BN}{{\mathbb N}}
\newcommand{\BP}{{\mathbb P}}
\newcommand{\BZ}{{\mathbb Z}}
\newcommand{\CA}{{\mathcal A}}

\newcommand{\CC}{{\mathcal C}}

\newcommand{\CF}{{\mathcal F}}
\newcommand{\CG}{{\mathcal G}}

\newcommand{\CI}{{\mathcal I}}
\newcommand{\CO}{{\mathcal O}}
\newcommand{\CP}{{\mathcal P}}

\newcommand{\CR}{{\mathcal R}}

\newcommand{\CX}{{\mathcal X}}
\newcommand{\CY}{{\mathcal Y}}
\newcommand{\CZ}{{\mathcal Z}}
\newcommand{\dd}{{\mathbf d}}

\newcommand{\ov}{\overline}

\newtheorem{thm}{Theorem}[section]
\newtheorem{cor}[thm]{Corollary}
\newtheorem{lem}[thm]{Lemma}
\newtheorem{prop}[thm]{Proposition}
\theoremstyle{definition}
\newtheorem{ex}[thm]{Example}
\numberwithin{equation}{section}
\begin{document}
\title{Codimension two singularities for representations of extended
 Dynkin quivers
 }
\author{Grzegorz Zwara}
%
\maketitle

\begin{abstract}
Let $M$ and $N$ be two representations of an extended Dynkin quiver 
such that the orbit $\CO_N$ of $N$ is contained in the orbit closure 
$\ov{\CO}_M$ and has codimension two.
We show that the pointed variety $(\ov{\CO}_M,N)$ is smoothly equivalent 
to a simple surface singularity of type $\BA_n$, or to the cone over 
a rational normal curve.   
\end{abstract}

\section{Introduction and the main results}

Throughout the paper, $k$ denotes an algebraically closed field, and 
$Q=(Q_0,Q_1,s,e)$ is a finite quiver, i.e.\ $Q_0$ is a finite set of 
vertices and $Q_1$ is a finite set of arrows $\alpha:s(\alpha)\to e(\alpha)$, 
where $s(\alpha)$ and $e(\alpha)$ denote the starting and the ending vertex 
of $\alpha$, respectively.
A representation $V$ of $Q$ over $k$ is a collection
$(V(i);\,i\in Q_0)$ of finite dimensional $k$-vector spaces
together with a collection
$(V(\alpha):V(s(\alpha))\to V(e(\alpha));\,\alpha\in Q_1)$ of
$k$-linear maps.
A morphism $f:V\to W$ between two representations is
a collection $(f(i):V(i)\to W(i);\,i\in Q_0)$ of $k$-linear maps such that
$$
f(e(\alpha))\circ V(\alpha)=W(\alpha)\circ f(s(\alpha))\qquad
\text{for all $\alpha\in Q_1$}.
$$
The dimension vector of a representation $V$ of $Q$ is the vector
$$
\dimv V=(\dim_k V(i))\in\BN^{Q_0}.
$$
We denote the category of representations of $Q$ by $\rep(Q)$,
and for any vector $\dd=(d_i)\in\BN^{Q_0}$
$$
\rep_Q(\dd)=\prod_{\alpha\in Q_1}
\BM_{d_{e(\alpha)}\times d_{s(\alpha)}}(k)
$$
is the vector space of representations $V$ of $Q$ with
$V(i)=k^{d_i}$, $i\in Q_0$.
The group
$$
\GL(\dd)=\prod_{i\in Q_0}\GL_{d_i}(k)
$$
acts on $\rep_Q(\dd)$ by
$$
((g_i)\star V)(\alpha)=g_{e(\alpha)}\cdot V(\alpha)\cdot 
g^{-1}_{s(\alpha)}.
$$
Given a representation $V$ of $Q$, we denote by $\CO_V$ the 
$\GL(\dd)$-orbit in $\rep_Q(\dd)$ consisting of the 
representations isomorphic to $V$, where $\dd=\dimv V$. 
An interesting problem is to study singularities of the
Zariski closure $\ov{\CO}_V$ of an orbit $\CO_V$ in $\rep_Q(\dd)$.

Following Hesselink (see~\cite[(1.7)]{Hes})
we call two pointed varieties $(\CX,x_0)$ and $(\CY,y_0)$
smoothly equivalent if there are smooth morphisms
$f:\CZ\to\CX$, $g:\CZ\to\CY$ and a point $z_0\in\CZ$
with $f(z_0)=x_0$ and $g(z_0)=y_0$.
This is an equivalence relation and the equivalence classes
will be denoted by $\Sing(\CX,x_0)$ and called the types
of singularities.
Obviously the regular points of the varieties form one type
of singularity, which we denote by $\Reg$.
Let $M$ and $N$ be representations in $\rep_Q(\dd)$ such that 
$M$ degenerates to $N$ ($N$ is a degeneration of $M$), i.e.\ 
$\CO_N\subseteq\ov{\CO}_M$.
We shall write $\Sing(M,N)$ for $\Sing(\ov{\CO}_M,n)$, where $n$
is an arbitrary point of $\CO_N$, and denote by $\codim(M,N)$ the
codimension of $\CO_N$ in $\ov{\CO}_M$.
We refer to \cite{BB}, \cite{Bmin}, \cite{Zsmo}, \cite{Zuni},
\cite{Zcodim1}, \cite{Zcodim2} and \cite{Zrat} for results 
in this direction.
Some of the results are expressed in terms of finite dimensional 
modules over finitely generated associative $k$-algebras, so it needs 
an explanation:
Given a representation $V$ of $Q$, we associate a (left) module 
$\widetilde{V}$ over the path algebra $kQ$ of $Q$, whose underlying 
vector space is $\bigoplus_{i\in Q_0}V(i)$.
This leads to an equivalence between $\rep(Q)$ and the category 
of finite dimensional $kQ$-modules.
Moreover, the equivalence preserves degenerations (of representations
and of modules, respectively) as well as their codimensions and 
types of singularities (see \cite{Bgeo}).
Applying \cite[Thm.1.1]{Zcodim1} (and the above geometric equivalence between 
representations of $Q$ and modules over $kQ$), we get $\Sing(M,N)=\Reg$ 
if $\codim(M,N)=1$.

We assume now that $\codim(M,N)=2$.
It was shown recently (\cite[Thm.1.3]{Zcodim2}) that $\Sing(M,N)=\Reg$ 
provided $Q$ is a Dynkin quiver.
This leads to a natural question about $\Sing(M,N)$ if $Q$ is an extended
Dynkin quiver, i.e.\ one of the following quivers
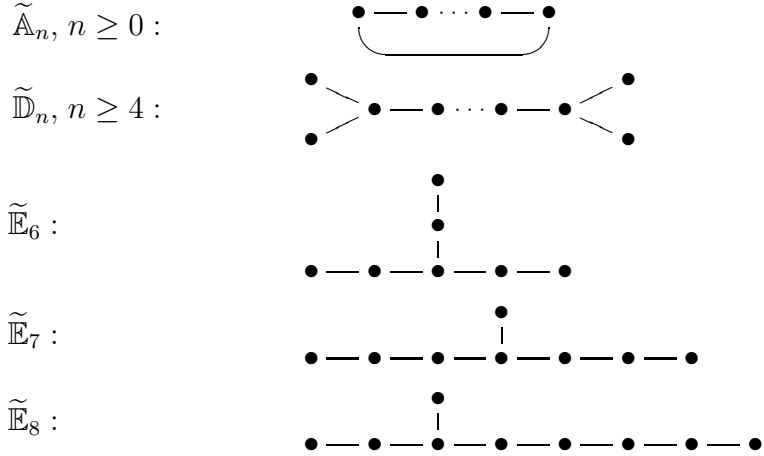
\begin{align*}
\widetilde{\BA}_n,&\;n\geq 0:&&\xymatrixcolsep{1pc}\xymatrixrowsep{.5pc}
 \vcenter{\xymatrix{
 &\bullet\ar@{-}`d[r]`[rrr][rrr]&\bullet\ar@{-}[l]&\bullet
 \ar@{}[l]|{\cdot\;\cdot\;\cdot}&\bullet\ar@{-}[l]
 }}\\
\widetilde{\BD}_n,&\;n\geq 4:&&\xymatrixcolsep{1pc}\xymatrixrowsep{0pc}
 \vcenter{\xymatrix{
 \bullet\ar@{-}[dr]&&&&&\bullet\\
 &\bullet\ar@{-}[r]&\bullet\ar@{}[r]|{\cdot\;\cdot\;\cdot}
 &\bullet\ar@{-}[r]&\bullet\ar@{-}[ur]\ar@{-}[dr]\\
 \bullet\ar@{-}[ur]&&&&&\bullet
 }}\\
\widetilde{\BE}_6:&&&\xymatrixcolsep{1pc}\xymatrixrowsep{.5pc}
 \vcenter{\xymatrix{
 &&\bullet\ar@{-}[d]\\
 &&\bullet\ar@{-}[d]\\
 \bullet\ar@{-}[r]&\bullet\ar@{-}[r]&\bullet\ar@{-}[r]&\bullet\ar@{-}[r]
 &\bullet
 }}\\
\widetilde{\BE}_7:&&&\xymatrixcolsep{1pc}\xymatrixrowsep{.5pc}
 \vcenter{\xymatrix{
 &&&\bullet\ar@{-}[d]\\
 \bullet\ar@{-}[r]&\bullet\ar@{-}[r]&\bullet\ar@{-}[r]&\bullet\ar@{-}[r]
 &\bullet\ar@{-}[r]&\bullet\ar@{-}[r]&\bullet
 }}\\
\widetilde{\BE}_8:&&&\xymatrixcolsep{1pc}\xymatrixrowsep{.5pc}
 \vcenter{\xymatrix{
 &&\bullet\ar@{-}[d]\\
 \bullet\ar@{-}[r]&\bullet\ar@{-}[r]&\bullet\ar@{-}[r]&\bullet\ar@{-}[r]
 &\bullet\ar@{-}[r]&\bullet\ar@{-}[r]&\bullet\ar@{-}[r]&\bullet
 }}
\end{align*}
(here $\xymatrixcolsep{1.5pc}\xymatrix@1{\bullet\ar@{-}[r]&\bullet}$ stands for
$\xymatrixcolsep{1.5pc}\xymatrix@1{\bullet\ar[r]&\bullet}$ or 
$\xymatrixcolsep{1.5pc}\xymatrix@1{\bullet&\bullet\ar[l]}$).
In the case of the Kronecker quiver
$$
Q=\quad\xymatrix@1{\bullet\ar@<.7ex>[r]\ar@<-.7ex>[r]&\bullet},
$$
two series $\BA_r=\Sing(\CA_{r+1},0)$, $\BC_r=\Sing(\CC_r,0)$, $r\geq 1$, 
of types of singularities occur (see \cite{BB}), where 
\begin{align*}
\CA_r&=\left\{(x,y,z)\in k^3;\;x^r=yz\right\}
 =\left\{(uv,u^r,v^r)\in k^3;\;u,v\in k\right\},\\
\CC_r&=\left\{(x_0,\ldots,x_r)\in k^{r+1};\;x_ix_j=x_lx_m\text{ if }
 i+j=l+m\right\}\\
&=\left\{(u^r,u^{r-1}v,\ldots,v^r)\in k^{r+1};\;u,v\in k\right\}.
\end{align*}
Thus $\BA_r$ is a simple surface singularity (a rational double point, 
a Kleinian or Du Val singularity), and $\CC_r$ is the affine cone over 
a rational normal curve of degree $r$.
Obviously $\BC_1=\Reg$, $\BC_2=\BA_1$ and the remaining types are pairwise 
different.
Note that, if $k$ is of characteristic zero, $\CA_r$ and $\CC_r$
are quotients of the plane $k^2$ by a cyclic subgroup of $\GL_2(k)$ isomorphic 
to $\BZ/r\BZ$.
We show that no other types of singularities can occur for representations 
of extended Dynkin quivers. 

\begin{thm} \label{mainthm}
Let $Q$ be an extended Dynkin quiver. 
Let $M$ and $N$ be representations in $\rep(Q)$ such that $M$
degenerates to $N$ and $\codim(M,N)=2$.
Then $\Sing(M,N)$ equals $\BA_r$ or $\BC_r$ for some $r\geq 1$.
\end{thm}

Among extended Dynkin quivers, cyclic quivers
$$
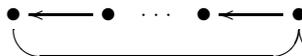
\xymatrix{
&\bullet\ar`d[r]`[rrr][rrr]&\bullet\ar[l]&\bullet
\ar@{}[l]|{\cdot\;\cdot\;\cdot}&\bullet\ar[l]
}
$$
play a special role.
For example, the path algebra $kQ$ is infinite dimensional and the category 
$\rep(Q)$ does not contain preprojective or preinjective representations.
For a basic background on the representation theory of extended Dynkin 
quivers we refer to \cite{Rin}.
We show that the types $\BC_r$, $r\geq 3$, do not occur
for cyclic quivers.

\begin{thm} \label{mainthm2}
Let $Q$ be a cyclic quiver. 
Let $M$ and $N$ be representations in $\rep(Q)$ such that $N$
is a degeneration of $M$ and $\codim(M,N)=2$.
Then $\Sing(M,N)$ equals $\Reg$ or $\BA_r$ for some $r\geq 1$.
\end{thm}

In order to prove the above theorems, we can apply reductions described
in \cite[Thm.1.1 and 1.2]{Zcodim2}.
Namely, we may assume that the representations $M$ and $N$ are disjoint
(i.e.\ they have no non-zero direct summands in common) and $\nu(N)\leq 2$,
where $\nu(V)$ is the number of summands in a decomposition of 
a representation $V$ as a direct sum of indecomposables.

We collect in Section~\ref{generalquiver} some fundamental properties 
of homomorphisms, extensions and degenerations of representations 
of quivers, and then we develop reductions for types of singularities 
following from \cite[(2.1)]{Bmin}.
Section~\ref{nilpotent} is devoted to the proof of Theorem~\ref{mainthm2}
in case the representations $M$ and $N$ are nilpotent. 
We recall in Section~\ref{proofs} some basic facts from representation theory
of extended Dynkin quivers and then we finish the proofs of our main results.

The author gratefully acknowledges support from the Polish
Scientific Grant KBN No.\ 1 P03A 018 27.

\section{Degenerations of quiver representations}
\label{generalquiver}

Let $V$ be a representation in $\rep_Q(\dd)$ for some $\dd\in\BN^{Q_0}$.
The isotropy group of $V$ can be identified with the group of automorphisms 
of $V$, and therefore
$$
\dim\CO_V=\dim\GL(\dd)-[V,V].
$$
Here and subsequently, 
$$
[V',V'']=\dim_k\Hom_Q(V',V'')\quad\text{and}\quad
{}^1[V',V'']=\dim_k\Ext^1_Q(V',V''),
$$ 
for any representations $V'$ and $V''$ of $Q$.
Consequently, 
\begin{equation} \label{codimformula}
\codim(M,N)=[N,N]-[M,M]
\end{equation}
for any representations $M$ and $N$ of $Q$ such that $M$ degenerates to $N$.
We shall need the following characterization of degenerations 
of representations (see \cite{Zgiv}).

\begin{prop} \label{characteriz}
Let $M$ and $N$ be representations of $Q$. Then $M$ degenerates to $N$
if and only if there is an exact sequence in $\rep(Q)$ of the form
$$
0\to Z\to Z\oplus M\to N\to 0
$$
for some representation $Z$.
Moreover, we may assume that $Z$ has a filtration 
$$
0=N_0\subset N_1\subset N_2\subset\cdots\subset N_h=Z
$$
with quotients $N_i/N_{i-1}$ isomorphic to $N$.
\end{prop}

As a direct consequence, we get well known facts that $M$ degenerates 
to $U\oplus V$ for any short exact sequence in $\rep(Q)$ of the form
$$
0\to U\to M\to V\to 0,
$$
and using the functors $\Hom_Q(-,Y)$ and $\Hom_Q(Y,-)$, that
\begin{equation} \label{standardineq}
[M,Y]\leq[N,Y]\qquad\text{and}\qquad[Y,M]\leq[Y,N]
\end{equation}
for any representation $Y$ of $Q$ (see for example \cite{Rie}).

By \cite{Bdeg} and \cite{Zext}, we get the following two propositions 
leading to better understanding of degenerations for extended Dynkin quivers.

\begin{prop}
Assume that $Q$ is an extended Dynkin quiver.
Let $M$ and $N$ be representations of $Q$ with $\dimv M=\dimv N$. 
Then the following conditions are equivalent:
\begin{enumerate}
\item[(1)] $M$ degenerates to $N$,
\item[(2)] $[M,Y]\leq[N,Y]$ for all $Y$ in $\rep(Q)$,
\item[(3)] $[Y,M]\leq[Y,N]$ for all $Y$ in $\rep(Q)$.
\end{enumerate}
\end{prop}

\begin{prop} \label{extfortame}
Assume that $Q$ is an extended Dynkin quiver.
If $N$ is a minimal degeneration of a representation $M$ (i.e.\ 
$\ov{\CO}_N\varsubsetneq\ov{\CO}_M$, but there is no representation 
$W$ with $\ov{\CO}_N\varsubsetneq\ov{\CO}_W\varsubsetneq\ov{\CO}_M$),
then there is an exact sequence in $\rep(Q)$ of the form
$$
0\to U\to M\to V\to 0
$$
with $N\simeq U\oplus V$. 
\end{prop}

\begin{lem} \label{seqcodim1}
Let $\sigma:\;0\to U\to M\to V\to 0$ be a short exact sequence 
in $\rep(Q)$ such that $\codim(M,N)=1$, where $N=U\oplus V$.
Then
$$
[U,M]=[U,N],\quad[M,V]=[N,V],\quad[Y,M]=[Y,N],\quad[M,Y]=[N,Y],
$$
for any direct summand $Y$ of $M$.
\end{lem}

\begin{proof}
Applying the functor $\Hom_Q(V,-)$ to the sequence $\sigma$, we get
$$
[V,N]=[V,U\oplus V]>[V,M].
$$ 
Assume that $M=Y\oplus Z$. It follows from \eqref{codimformula}
that
\begin{align*}
1=([N,N]-[M,M])=&([U,N]-[U,M])+([V,N]-[V,M])\\
&+([N,Y]-[M,Y])+([N,Z]-[M,Z]).
\end{align*}
Thus $[U,N]=[U,M]$ and $[N,Y]=[M,Y]$.
In much the same way one can show that $[N,V]=[M,V]$ and 
$[Y,N]=[Y,M]$.
\end{proof}

We shall need the following sufficient condition for the regularity
of points in orbit closures.

\begin{cor} \label{suffreg}
Let $\sigma:\;0\to Y\oplus U\to Y\oplus M\to V\to 0$ be a short 
exact sequence in $\rep(Q)$ such that 
$[Y\oplus U\oplus M,M]=[Y\oplus U\oplus M,U\oplus V]$.
Then $\Sing(M,U\oplus V)=\Reg$.
\end{cor}

\begin{proof}
We apply \cite[Prop.2.2]{Zuni} to the direct sum of $\sigma$ and the short
exact sequence $0\to 0\to U\xrightarrow{\sim}U\to 0$.
\end{proof}

Let $0\to U\to M\to V\to 0$ be a short exact sequence in $\rep(Q)$.
We say that $V$ is a generic quotient of $M$ by $U$ if the orbit $\CO_V$
is dense in the set of representations in $\rep_Q(\dimv V)$ isomorphic
to the cokernels of monomorphisms from $U$ to $M$.
We say that $M$ is a generic extension of $V$ by $U$ if the orbit $\CO_M$
is dense in the set of representations in $\rep_Q(\dimv M)$ isomorphic
to the extensions of $V$ by $U$.  
Note that the above sets of all possible cokernels or extensions 
are constructible and irreducible (see \cite[(2.1)]{Bmin}).
We shall need the following modification of \cite[Thm.2.1]{Bmin}.

\begin{prop} \label{citedthm21}
Let $U$, $M$, $M'$, $V$ and $V'$ be representations of $Q$ satisfying
the following conditions:
\begin{enumerate}
\item[(1)] $M$ degenerates to $M'$,
\item[(2)] $[U,M]=[U,M']$,
\item[(3)] $V$ and $V'$ are the generic quotients of $M$ and $M'$, 
 respectively, by $U$,
\item[(4)] $M$ is the generic extension of $V$ by $U$.
\end{enumerate}
Then $V$ degenerates to $V'$, $\codim(V,V')\leq\codim(M,M')$ and 
$$
\Sing(M,M')=\Sing(V,V').
$$
\end{prop}

\begin{proof}
One can repeat the proof of \cite[Thm.2.1]{Bmin} with two differences.
First we omit an assumption that
$$
{}^1[V,U]-[V,U]={}^1[V',U]-[V',U].
$$
This equation was used in the proof of the cited theorem only to conclude 
that some map $p'$ was a vector bundle, but this map in the case of 
representations of a quiver (instead of modules over an algebra)
is in fact a trivial vector bundle.
The second difference is that we assume in addition that the quotient $V'$ 
of $M'$ by $U$ is generic.
Reading carefully the proof, we see that our additional assumption 
implies $\codim(V,V')\leq\codim(M,M')$.
\end{proof}

Let $M$ and $N$ be representations of $Q$ such that $M$ degenerates 
to $N$.
We want to apply the above proposition for $U$ being the socle 
$\soc(M)$ of $M$.
We note that $U$ is isomorphic to a direct summand of $\soc(N)$.
Indeed, the multiplicity of a simple representation $S$ of $Q$ 
as a direct summand of $\soc(V)$ equals $[S,V]$, 
for any representation $V$ of $Q$; 
and $[S,M]\leq[S,N]$, by \eqref{standardineq}.
Thus $\soc(N)\simeq U\oplus W$ for some representation $W$ of $Q$.
If the semi-simple representations $U$ and $W$ are disjoint, then 
there is a unique subrepresentation $U'$ of $N$ isomorphic to $U$,
as $U'$ must be contained in $\soc(N)$. 
In such a case we write $N/U$ for the quotient of $N$ by $U'$. 
  
\begin{cor} \label{corgenquiver}
Let $M$ and $N$ be representations of $Q$ such that $N$ is 
a degeneration of $M$.
Assume that $U=\soc(M)$ and its direct complement in $\soc(N)$ are 
disjoint representations.
Let $V=M/U$ and $V'=N/U$.
If $M$ is the generic extension of $V$ by $U$ then:
\begin{enumerate}
\item[(1)] $V$ degenerates to $V'$,
\item[(2)] $\codim(V,V')\leq\codim(M,N)$, 
\item[(3)] $\Sing(V,V')=\Sing(M,N)$.
\end{enumerate}
\end{cor}

\begin{proof}
Since $U$ is a semisimple representation, we get
$$
[U,M]=[U,\soc(M)]=[U,\soc(N)]=[U,N].
$$ 
Obviously $V$ is the generic quotient of $M$ by $U$, 
and $V'$ is the generic quotient of $N$ by $U$.
Thus the claim follows from Proposition~\ref{citedthm21}.
\end{proof}

\section{Nilpotent representations of cyclic quivers}
\label{nilpotent}

We fix a positive integer $n$.
Let $Q$ be the cyclic quiver with the set of vertices $Q_0=\BZ/n\BZ$ 
and the set of arrows 
$Q_1=\{\alpha_{\ov{l}}:\ov{l}\to\ov{l-1};\;\ov{l}\in\BZ/n\BZ\}$:
$$
\xymatrix@1{
\ov{1}\ar `d[r] `[rrr]_{\alpha_{\ov{1}}}[rrr]&\ov{2}\ar[l]_{\alpha_{\ov{2}}}
&\cdots\ar[l]_{\alpha_{\ov{3}}}&\ov{n}\ar[l]_{\alpha_{\ov{n}}}
}
$$
We call a representation $V=(V(\ov{l}),V(\alpha_{\ov{l}}))_{\ov{l}\in\BZ/n\BZ}$ 
of $Q$ \textit{nilpotent} if the endomorphisms
$$
V(\alpha_{\ov{l-n+1}})\circ\cdots\circ V(\alpha_{\ov{l-1}})\circ
V(\alpha_{\ov{l}}):V(\ov{l})\to V(\ov{l}),
\quad\ov{l}\in\BZ/n\BZ,
$$
are nilpotent, or equivalently, if there is a positive integer $h$ such that
$$
V(\alpha_{\ov{l-h+1}})\circ\cdots\circ V(\alpha_{\ov{l-1}})\circ
V(\alpha_{\ov{l}})=0
$$
for any $\ov{l}\in\BZ/n\BZ$.
We denote by $\rep^0(Q)$ the full subcategory of $\rep(Q)$ 
of nilpotent representations.
It is an abelian subcategory closed under extensions.
The aim of the section is to prove the following result.

\begin{prop} \label{keynilpotent}
Let $M$ and $N$ be nilpotent representations of $Q$ such that $N$ is 
a degeneration of $M$ and $\codim(M,N)=2$.
Then $\Sing(M,N)$ equals $\Reg$ or $\BA_r$ for some $r\geq 1$.
\end{prop}

For any two integers $i\leq j$ we consider an indecomposable nilpotent 
representation $V_{i,j}$ described by a basis 
$$
\{b_i,b_{i+1},\ldots,b_j\}\subset\bigoplus_{\ov{l}\in\BZ/n\BZ}
 V_{i,j}(\ov{l}),\quad
b_l\in V_{i,j}(\ov{l}),\;
V_{i,j}(\alpha_{\ov{l}})(b_l)=\begin{cases}b_{l-1}&l>i,\\ 0&l=i.\end{cases}
$$
Observe that $\dim_kV_{i,j}=j-i+1$.
Any indecomposable nilpotent representation of $Q$ is isomorphic to some 
$V_{i,j}$, and $V_{i,j}$ is isomorphic to another $V_{i',j'}$ 
if and only if $i'=i+cn$ and $j'=j+cn$ for some integer $c$.

Observe that $S_i=V_{i,i}$ is a simple representation of $Q$ supported 
at the vertex $\ov{i}$ for any $i\in\BZ$.
Moreover, 
\begin{equation} \label{socfornilpotent}
\soc(V_{i,j})\simeq S_i\qquad\text{and}\qquad
V_{i,j}/\soc(V_{i,j})\simeq\begin{cases}0&i=j,\\ V_{i+1,j}&i<j,
\end{cases}
\end{equation}
for any integers $i\leq j$.

\begin{lem} \label{uniquelybysoc}
Let $V$ and $W$ be two nilpotent representations of $Q$ such that
$\soc(V)\simeq\soc(W)$ and $V/\soc(V)\simeq W/\soc(W)$.
Then $V\simeq W$.
\end{lem}

\begin{proof}
Let $v_{i,j}$, $t_{i,j}$ and $u_{i,j}$ denote the multiplicities of 
$V_{i,j}$ as direct summands of $V$, $V/\soc(V)$ and $\soc(V)$, respectively.
It suffices to show that the numbers $v_{i,j}$'s depend only on $t_{i,j}$'s
and $u_{i,j}$'s.
By \eqref{socfornilpotent}, we get $v_{i,j}=t_{i+1,j}$ provided $i<j$,
and $u_{i,i}=\sum_{j\geq i} v_{i,j}$.
Consequently, $v_{i,i}=u_{i,i}-\sum_{j>i}t_{i+1,j}$, and the claim follows.
\end{proof}

\begin{lem} \label{gensocnilpotent}
Let $M$ be a nilpotent representation of $Q$. Then $M$ is a generic
extension of $M/\soc(M)$ by $\soc(M)$.
\end{lem}

\begin{proof}
The category $\rep^0(Q)$ is closed under extension, hence there is
up to isomorphism only finitely many extensions of $M/\soc(M)$ by $\soc(M)$.
Since the set of representations in $\rep_Q(\dimv M)$ isomorphic to
extensions of $M/\soc(M)$ by $\soc(M)$ is irreducible, there exists the 
generic extension $E$. 
In particular, $E$ degenerates to $M$ and $\soc(E)$ is isomorphic to 
a direct summand of $\soc(M)$ (see Section~\ref{generalquiver}).
On the other hand, we conclude from the short exact sequence
$$
0\to\soc(M)\to E\to M/\soc(M)\to 0
$$
that $\soc(M)$ is isomorphic to a subrepresentation of $\soc(E)$.
Hence $\soc(E)$ is isomorphic to $\soc(M)$ and $E/\soc(E)$ is isomorphic to
$M/\soc(M)$.
Consequently, $E$ is isomorphic to $M$, by Lemma~\ref{uniquelybysoc}.
\end{proof}

We say that a pair $(M,N)$ of nilpotent representations of $Q$
is \textit{admissible} if $M$ degenerates to $N$, $\codim(M,N)\leq 2$ 
and $\nu(N)\leq 2$.
Combining Corollary~\ref{corgenquiver}, Lemma~\ref{gensocnilpotent}
and the fact that $\nu(V/\soc(V))\leq\nu(V)$ for any nilpotent 
representation $V$ of $Q$, we get the following result.

\begin{cor} \label{reduct1}
Let $(M,N)$ be an admissible pair of nonzero nilpotent representations 
of $Q$.
Then one of the following conditions holds:
\begin{enumerate}
\item[(1)] $\soc(M)\simeq S_i$ and $\soc(N)\simeq S_i\oplus S_i$ 
for some integer $i$,
\item[(2)] there is an admissible pair $(M',N')$ with 
$\dim_kM'<\dim_kM$ and 
$$
\Sing(M',N')=\Sing(M,N).
$$
\end{enumerate}
\end{cor}

Now we consider the radical and the top of nilpotent representations.
Observe that
$$
\rad(V_{i,j})\simeq\begin{cases}0&i=j,\\ V_{i,j-1}&i<j,\end{cases}
\qquad\text{and}\qquad\top(V_{i,j})=V_{i,j}/\rad(V_{i,j})\simeq S_j,
$$
for any integers $i\leq j$.
By duality, we obtain the following result.

\begin{cor} \label{reduct2}
Let $(M,N)$ be an admissible pair of nonzero nilpotent representations 
of $Q$.
Then one of the following conditions holds:
\begin{enumerate}
\item[(1')] $\top(M)\simeq S_j$ and $\top(N)\simeq S_j\oplus S_j$ 
for some integer $j$,
\item[(2)] there is an admissible pair $(M',N')$ with 
$\dim_kM'<\dim_kM$ and 
$$
\Sing(M',N')=\Sing(M,N).
$$
\end{enumerate}
\end{cor}

\begin{proof}[Proof of Proposition~\ref{keynilpotent}]
By \cite[Thm.1.1 and 1.2]{Zcodim2}, we may assume that $\nu(N)\leq 2$
(as mentioned in Section 1), which implies that the pair $(M,N)$ 
is admissible.
Applying the reductions described in Corollaries~\ref{reduct1} 
and~\ref{reduct2} as many times as possible, we may assume that 
the conditions (1) and (1') hold (otherwise $\Sing(M,N)=\Sing(0,0)=\Reg$).
Thus, $M\simeq V_{i,j+an}$ and $N\simeq V_{i,j+bn}\oplus V_{i,j+cn}$
for some integers $i$, $j$, $a$, $b$ and $c$.
Without loss of generality we may assume that $i-n\leq j<i$ 
as $S_l=S_{l+n}$ for any integer $l$.
We conclude from the equalities
$$
j-i+an+1=\dim_kM=\dim_kN=(j-i+bn+1)+(j-i+cn+1)
$$
that $j=(i-1)+(a-b-c)n$. 
Hence $j=i-1$, $a=b+c$ and the numbers $a$, $b$ and $c$ are positive.

Let $f$ be a positive integer and $Q'$ be a loop quiver with a unique 
arrow $\gamma$.
Let $U_f$ denote the representation in $\rep_{Q'}(f)$ such that 
$U_f(\gamma)$ is the nilpotent Jordan block matrix (of size $f$).   
Observe that up to isomorphism,
$V_{i,j+fn}(\alpha_{\ov{i}})$ is the nilpotent Jordan block matrix 
of size $f$ and $V_{i,j+fn}(\beta)$ is the identity matrix of size $f$
for the remaining arrows $\beta$ in $Q_1$.
Hence using the operation ``replacing one arrow by none'', described 
in \cite[(5.2)]{Bmin}, to the arrows $\beta\neq\alpha_{\ov{i}}$, we 
conclude that $\codim(M,N)=\codim(U_a,U_b\oplus U_c)$ and
$$
\Sing(M,N)=\Sing(U_a,U_b\oplus U_c).
$$
Observe that $[U_f,U_g]=\min\{f,g\}$ for any positive integers $f$ 
and $g$.
Thus
\begin{align*}
2&\geq\codim(U_a,U_b\oplus U_c)=[U_b\oplus U_c,U_b\oplus U_c]
 -[U_a,U_a]\\
&=b+c+2\min\{b,c\}-a=2\min\{b,c\},
\end{align*}
which implies that $\min\{b,c\}=1$.
We may assume that $b=1$. Hence the claim follows from a well known 
fact that $\Sing(U_{c+1},U_1\oplus U_c)=\BA_c$ 
(for instance, see \cite{KP} or \cite[(2.2)]{Bmin}).
\end{proof}

\begin{ex}
We shall illustrate the reductions used in the proof of 
Proposition~\ref{keynilpotent} for $n=2$. 
Let $M=V_{1,4}$ and $N=V_{1,2}\oplus V_{2,3}$.
One can show that $M$ degenerates to $N$ and $\codim(M,N)=2$.
Using the first reduction and then three times the second one we get
\begin{align*}
\Sing(V_{1,4},V_{1,2}\oplus V_{2,3})
 &=\Sing(V_{2,4},V_{2,2}\oplus V_{2,3})=\Sing(V_{2,3},V_{2,3})\\
&=\Sing(V_{2,2},V_{2,2})=\Sing(0,0)=\Reg.
\end{align*}
It is not difficult to see that $\codim(V_{2,4},V_{2,2}\oplus V_{2,3})=1$.

Now let $M=V_{1,1}\oplus V_{2,8}$ and $N=V_{1,3}\oplus V_{2,6}$.
Then $M$ degenerates to $N$, $\codim(M,N)=2$ and
\begin{align*}
\Sing(V_{1,1}\oplus V_{2,8}, V_{1,3}\oplus V_{2,6})
 &=\Sing(V_{3,8},V_{2,3}\oplus V_{3,6})
 =\Sing(V_{4,8},V_{2,3}\oplus V_{4,6})\\
&=\Sing(V_{4,7},V_{2,3}\oplus V_{4,5})=\Sing(V_{0,3},V_{0,1}\oplus V_{0,1})\\
&=\Sing(U_2,U_1\oplus U_1)=\BA_1.
\end{align*}
\end{ex}

We shall need a fact that geometric properties of orbit closures are
preserved if we pass from $\rep^0(Q)$ to an equivalent exact category. 

\begin{prop} \label{nilptostabletube}
Let $\CF:\rep^0(Q)\to\CA$ be an equivalence of exact subcategories,
where $\CA$ is a full subcategory  closed under extensions of $\rep(Q')$
for some quiver $Q'$.
Let $M$ and $N$ be two representations in $\rep^0(Q)$.
Then $M$ degenerates to $N$ if and only if $\CF(M)$ degenerates 
to $\CF(N)$. 
Moreover, if this is the case, then $\codim(\CF(M),\CF(N))=\codim(M,N)$
and
$$
\Sing(\CF(M),\CF(N))=\Sing(M,N).
$$
\end{prop}

\begin{proof}
The first part follows from Proposition~\ref{characteriz},
as the equivalence $\CF$ is an exact functor and the subcategories 
$\rep^0(Q)$ and $\CA$ are closed under extensions.
Thus we assume that $M$ degenerates to $N$ and $\CF(M)$ degenerates
to $\CF(N)$.
The equality of codimensions follows from \eqref{codimformula}.
Let $\rep^{0,h}(Q)$ denote the full subcategory of $\rep^0(Q)$ 
consisting of the representations $V$ such that $V(\omega)=0$
for any  path in $Q$ of length $h\geq 1$.
We choose $h$ such that $M$ and $N$ belong to $\rep^{0,h}(Q)$
(for example, $h=\dim_kM=\dim_kN$).
Let $\CG:\rep^{0,h}(Q)\to\rep(Q')$ be a restriction of $\CF$
followed by the inclusion of $\CA$ in $\rep(Q')$.
The category $\rep^{0,h}(Q)$ is equivalent to the category of
modules over some finite dimensional algebra $B$ and the functor
$\CG$ is hom-controlled in the sense of \cite{Zsmo}.
Hence 
$$
\Sing(\CF(M),\CF(N))=\Sing(\CG(M),\CG(N))=\Sing(M,N),
$$
by \cite[Thm.1.2]{Zsmo} and the geometric equivalence (\cite{Bgeo})
between representations in $\rep^{0,h}(Q)$ and $B$-modules.
\end{proof}

\section{Proof of the main results}
\label{proofs}

Throughout the section, $Q$ is an extended Dynkin quiver, and $M$, $N$ are
representations of $Q$ such that $M$ degenerates to $N$ and $\codim(M,N)=2$. 
In order to prove the theorems, we may assume that the representations
$M$ and $N$ are disjoint and $\nu(N)\leq 2$.
Let $W$ be a degeneration of $M$ such that $N$ is a minimal degeneration 
of $W$.
It follows from Proposition~\ref{extfortame} that there is a short exact
sequence 
$$
\sigma:\quad 0\to U\to W\to V\to 0
$$
in $\rep(Q)$ such that $N\simeq U\oplus V$.
Thus the above sequence does not split, ${}^1[V,U]>0$,
$\nu(N)=2$, and the representations $U$ and $V$ are indecomposable. 
Moreover, applying \eqref{standardineq}, and the functors $\Hom_Q(-,U)$ 
and $\Hom_Q(V,-)$ to $\sigma$ we get
\begin{equation} \label{UleftVright}
[N,U]>[W,U]\geq[M,U]\quad\text{and}\quad[V,N]>[V,W]\geq[V,M].
\end{equation}

We need to recall a few facts and definitions from \cite[(3.6)]{Rin}.
Assume first that $Q$ is not a cyclic quiver, or equivalently, $Q$ has
no oriented cycles.
The category $\rep(Q)$ decomposes into three exact subcategories 
$\CP$, $\CI$ and $\CR$, consisting of the preprojective, preinjective 
and regular representations, respectively.
The category $\CR$ is abelian and decomposes further into a $\BP^1(k)$-family
$\coprod_{\lambda\in\BP^1(k)}\CR_\lambda$ of uniserial categories.
The category $\CR_\lambda$ is equivalent to the category of nilpotent
representations of a cyclic quiver with $r_\lambda\geq 1$ vertices,
considered already in the previous section.
Now assume that $Q$ is a cyclic quiver.
Then the description of the category $\rep(Q)$ is even simpler.
Namely, $\CP=\CI=0$ and $\rep(Q)=\CR$ decomposes into a $k$-family 
$\coprod_{\lambda\in k}\CR_\lambda$, where $\CR_0$ consists of the nilpotent 
representations, and $\CR_\lambda$, for $\lambda\neq 0$,
is equivalent to the category of nilpotent representations of a loop 
quiver ($r_\lambda=1$).
The following lemma contains important information on homomorphisms 
and extensions for representations of $Q$.
 
\begin{lem} \label{wherearehomext}
Assume that $X$ and $Y$ are indecomposable representations of $Q$, such that
$[X,Y]>0$ or ${}^1[Y,X]>0$.
Then $X$ is preprojective, or $Y$ is preinjective, or both representations
belong to $\CR_\lambda$ for some $\lambda\in\BP^1(k)$.
\end{lem}

The following corollary finishes the proof of Theorem~\ref{mainthm2}.

\begin{cor}
If the representation $N$ is regular then $\Sing(M,N)$ equals $\Reg$ 
or $\BA_r$ for some $r\geq 1$.
\end{cor}

\begin{proof}
Since ${}^1[V,U]>0$, both representations belong to some 
$\CR_\lambda$.
Let $Y$ be an indecomposable direct summand of $M$.
Using \eqref{standardineq}, we get
$$
[U\oplus V,Y]\geq[M,Y]>0
\qquad\text{and}\qquad
[Y,U\oplus V]\geq[Y,M]>0.
$$
Hence $Y$ must belong to $\CR_\lambda$, by Lemma~\ref{wherearehomext}.
This implies that $M\oplus N$ belongs to the category $\CR_\lambda$,
and the claim follows from Propositions~\ref{keynilpotent}
and~\ref{nilptostabletube}.
\end{proof}

From now on, we assume that the quiver $Q$ is not cyclic, and $N$ has 
a nonzero preprojective direct summand (the case $N$ has a nonzero 
preinjective direct summand follows by duality).  
Let $\ind(\CP)$ denote a complete set of pairwise non-isomorphic 
indecomposable preprojective representations of $Q$.
There is a partial order $\preceq$ on $\ind(\CP)$ such that 
$[X,Y]>0$ implies $X\preceq Y$ for any $X$ and $Y$ in $\ind(\CP)$.
By \cite[Lem.3.1]{Bext}, there is a $\preceq$-minimal 
$T$ in $\ind(\CP)$ with the property $[N,T]>[M,T]$, and any such $T$ 
is a direct summand of $N$.
Moreover, using the Auslander-Reiten formula mentioned in the
proof of \cite[Lem.3.1]{Bext}, we conclude that $[T,N]=[T,M]$.
By \eqref{UleftVright}, $T$ is not isomorphic to $V$.
Thus $T\simeq U$ and
\begin{equation} \label{fromU}
[U,N]=[U,M].
\end{equation}
If $[N,V]=[M,V]$, then $\Sing(M,N)=\BC_r$ for some $r\geq 1$, 
by \cite[Thm.1.1]{Zrat}.
Hence we may assume that
\begin{equation} \label{toV}
[N,V]>[M,V].
\end{equation}
We shall show that $\Sing(M,N)=\Reg$.
By \eqref{codimformula},
$$
2=([N,U]-[M,U])+([N,V]-[M,V])+([M,N]-[M,M]).
$$
Combining this equality with \eqref{UleftVright} and \eqref{toV}, we get
$$
[N,U]-[M,U]=[N,V]-[M,V]=1\quad\text{and}\quad[M,N]=[M,M].
$$
Using the equality \eqref{fromU} gives
\begin{equation} \label{fromUM}
[U\oplus M,N]=[U\oplus M,M].
\end{equation}
If $M\simeq W$, then $\Sing(M,N)=\Reg$, by Corollary~\ref{suffreg}
applied to $\sigma$.

From now on, we assume that $W$ is not isomorphic to $M$, i.e.\ $W$
is a proper degeneration of $M$.
Then $\codim(M,W)=\codim(W,N)=1$.
In particular, $W$ is a minimal degeneration of $M$, and there
is a short exact sequence
$$
\eta:\quad 0\to W'\to M\to W''\to 0
$$
in $\rep(Q)$ with $W'\oplus W''\simeq W$, 
by Proposition~\ref{extfortame}.
Applying Lemma~\ref{seqcodim1} to the exact sequences $\eta$ 
and $\sigma$, we get 
\begin{equation} \label{fromW'}
[W',M]=[W',W]=[W',N].
\end{equation}
Considering the sequence $\sigma$ and the direct sum of the sequence $\eta$
and 
$$
0\to 0\to W'\xrightarrow{\sim}W'\to 0,
$$
we get the following commutative diagram with exact rows and columns
$$
\xymatrix{
&&0\ar[d]&0\ar[d]\\
&&W'\ar@{=}[r]\ar[d]&W'\ar[d]\\
\theta:&0\ar[r]&X\ar[r]\ar[d]&W'\oplus M\ar[r]\ar[d]
 &V\ar[r]\ar@{=}[d]&0\\
\sigma:&0\ar[r]&U\ar[r]\ar[d]&W'\oplus W''\ar[r]\ar[d]&V\ar[r]&0\\
&&0&0
}
$$
for some representation $X$.
Applying the functor $\Hom_Q(U,-)$ to $\theta$ and to the short 
exact sequence
$$ 
\psi:\quad 0\to W'\to X\to U\to 0,
$$
we obtain two non-negative integers
$$
[U,X\oplus V]-[U,W'\oplus M]\quad\text{and}\quad
[U,W'\oplus U]-[U,X].
$$
These numbers are zero as their sum equals
$$
[U,V]-[U,M]+[U,U]=[U,N]-[U,M]=0,
$$
by \eqref{fromU}.
Consequently, the last map in the exact sequence
$$
0\to\Hom_Q(U,W')\to\Hom_Q(U,X)\to\Hom_Q(U,U)
$$ 
induced by $\psi$ is surjective.
This implies that the exact sequence $\psi$ splits.
Thus $X$ is isomorphic to $W'\oplus U$.
Using the equalities \eqref{fromUM} and \eqref{fromW'}, we get
$\Sing(M,N)=\Reg$, by Proposition~\ref{suffreg} applied to 
the sequence $\psi$.
This finishes the proof of Theorem~\ref{mainthm}.
\qed


\bigskip

\noindent
Grzegorz Zwara\\
Faculty of Mathematics and Computer Science\\
Nicolaus Copernicus University\\
Chopina 12/18, 87-100 Toru\'n, Poland\\
E-mail: gzwara@mat.uni.torun.pl
\end{document}